\documentclass{article}
\usepackage{authblk}
\usepackage[british]{babel}
\usepackage{amsmath,graphicx,amssymb, amsfonts, amsthm}
\usepackage{subfigure}
\usepackage[utf8]{inputenc}

\title{Uncertainty quantification of modified Cahn-Hilliard equation for image inpainting}
\author{Yin Xian}
\affil{Department of Mathematics, Hong Kong University of Science and Technology, Clear Water Bay, Hong Kong}
\date{June 2019}

\begin{document}

\maketitle
\section*{Abstract}
In this paper, we review modified Cahn-Hilliard equation for image inpainting and explore the effect when the initial condition is uncertain. We study the statistical properties of the solution when the noise is present. The generalized polynomial chaos and the perturbation expansion are used to analyze the equation. Experimental results are attached for comparison of solution behavior. 

\section{Introduction}
\subsection{PDE based image inpainting}
Given an observed image $g$, which is corrupted, we want to have the original image $u$.
\begin{align*}
Tu=g
\end{align*}
where $T$ models the process through which the image $u$ went through observation. Let $\Omega \in \mathbb{R}^2$  be a given domain, $B_1$ be banach spaces over $\Omega$ and $g\in B_1$ be the given image.  A general variational approach in image processing can be written as:
\begin{align}
J(u)=R(u)+\frac{1}{2\lambda}||Tu-g ||^2_{B_1}
\end{align}
where $\lambda$ is the tuning parameter of the problem and $T\in \mathcal{L}(B_1)$ is a bounded linear operator.  $R$ denotes the regularizing term which smoothes the image $u$ and represents some kind of a priori information about the minimizer $u$. $||Tu-g||_{B_1}^2$ is called fidelity term of the approach which forces the minimizer $u$ to stay close to the given image $g$.

For $B_1=L^2(\Omega)$ we also have the corresponding Euler-Lagrange equation:
\begin{align}
-\lambda\nabla R(u)+T^{*}(g-Tu)=0,~~~~\text{in~} \Omega
\end{align}
the corresponding steepest descent equation for $u$ is the given image
\begin{align*}
u_t=-\lambda\nabla R(u)+T^*(g-Tu),~~~\text{in~} \Omega
\end{align*}
For Cahn-Hilliard equation and TV-$H^{-1}$ inpainting, the image processing approach is directly given by an evolutionary PDE. A regularizing term which can delete noise and preserve important image features like edges is the total variation. 

\begin{align}
J(u)=R(u)+\frac{1}{2\lambda}||\chi_{\Omega\backslash D}(u-g)||_{B_1}^{2}
\end{align}

where 
\begin{align*}
\chi_{\Omega\backslash D}(x)=\left\{
                \begin{array}{ll}
                  1, ~~ \Omega \backslash D  \\
                  0, ~~ \in D
                \end{array}
              \right.
\end{align*}
$R(u)$ fills in the image content into the missing domain $D$, by diffusion and transport. The fidelity term only has impact on the minimizer $u$ outside of the inpainting domain due to the characteristic function $\chi_{\Omega\backslash D}$. For $R(u)$ we have,

\begin{itemize}  
\item
$R(u)=\int_{\Omega}|\nabla u|^2dx$, harmonic inpainting. 
\item
$R(u)=\int_{\Omega}|\nabla u|dx$, TV-inpainting, (Chan and Shen 2001).
\item
$R(u)=\int_{\Omega}(1+\nabla\cdot(\frac{\nabla u}{|\nabla u|})))|\nabla u|dx$, Euler's elastica inpainting
\item
inpainting for binary images with the Cahn-Hilliard equation, (Bertozzi, Esedoglu and Gillette 06)
\item
TV-$H^{-1}$ inpainting, (Burger, He).
\end{itemize}  

For TV inpainting, it propagate sharp edges into the damaged domain
\begin{align}
\min\limits_{u}\int_{\Omega}|\nabla u|dx \Longleftrightarrow \min\limits_{\Gamma_{\lambda}}\int_{-\infty}^{\infty} \text{length}(\Gamma_{\lambda})d\lambda
\end{align}
where $\Gamma_{\lambda}=\{x\in\Omega: u(x)=\lambda\}$. It penalizes length of edges, cannot connect contours across very large distances, it can result in corners of the level lines across the inpainting domain.

For higher order approaches,  often they do not posses a maximum principle or comparison principle. For the proof of well-posedness of higher order inpainting models variational methods are often not applicable, and it need stable and fast numerical solvers. 
\begin{align*}
\min\limits_{u}\int_{\Omega}\left(a+b\nabla\cdot\left(\frac{\nabla u}{|\nabla u|}\right)\right) |\nabla u|dx \Longleftrightarrow \min\limits_{\Gamma_{\lambda}}\int_{-\infty}^{\infty} (a~\text{length}(\Gamma_{\lambda})+b~\text{curvature}(\Gamma_{\lambda}))d\lambda
\end{align*}

\subsection{Modified Cahn-Hilliard equation}
% bosch PhD thesis
(Analytical challenges) Results for stationary solution are difficult because of the missing energy for the equation.
\begin{align}
u_t=\Delta(-\varepsilon\Delta u+\frac{1}{\varepsilon}W'(u))+\lambda\chi_{\Omega\backslash D}(f-u)
\end{align}
$W(u)$ is a double well potential. The two wells of $W$ correspond to values of $u$ that are taken by most of the grey scale values. $\lambda>>0$. By choosing a potential with wells at the values 0 (black) and 1 (white), therefore provides a simple model for the inpainting of binary images. The parameter $\varepsilon>0$ determines the steepness of the transition between 0 and 1. 

The Cahn-Hilliard equation is a relatively simple fourth order PDE used for this task rather than more complex models involving curvature terms such as Euler-Elastica inpainting. It has many of the desirable properties of curvature-based inpainting models such as the smooth continuation of level lines into the missing domain~\cite{schonlieb}. 

The gradient flow in $H^{-1}$ of the energy is given by
\begin{align*}
J^{1}(u)=\int_{\Omega}\frac{\varepsilon}{2}|\nabla u|^2+\frac{1}{\varepsilon}W(u) dx.
\end{align*}
The gradient flow in $L^2$ of the energy is given by
\begin{align*}
J^2(u)=\frac{\lambda}{2}\int_{\Omega}\chi_{\Omega\backslash D}(f-u)^2dx
\end{align*}

The global existence for the evolution equation is given by Bertozzi, Esedoglu and Gillette 06, the authors proved that in the limit $\lambda\to \infty$, a stationary solution solves
\begin{align*}
\nabla(\varepsilon\nabla u-\frac{1}{\varepsilon}W'(u))&=0,~~\text{ in}~~D, \\
u&=f~~\text{ on}~~\partial D,  \\
\nabla u&=\nabla f~~\text{ on}~~\partial D,
\end{align*} 
for $f$ regular enough ($f\in C^2$). The existence of a stationary solution to the modified Cahn-Hilliard equation is proved by using the idea of fixed point equation.  This claim that fourth-order methods are superior to second order methods with respect to a smooth continuation of the image contents into the missing domain. Combined with other inpainting strategies, it can constitute a powerful method for inpainting of the structural part of an image~\cite{schonlieb}. 

The sequence of Cahn-Hilliard functionals:
\begin{align*}
\int_{\Omega}(\frac{\varepsilon}{2}|\nabla u|^2+\frac{1}{\varepsilon}W(u))dx
\end{align*}
$\Gamma$-convergence in the topology $L^1(\Omega)$ to
\begin{align}
TV(u)=\left\{
                \begin{array}{ll}
                  C_0|Du|(\Omega), ~~ \text{if $|u(x)|=1$ in $\Omega$}  \\
                  +\infty, ~~ otherwise
                \end{array}
              \right.
\end{align}
as $\varepsilon\to 0$, where $C_0=2\int_{-1}^{1}\sqrt{F(s)}ds$. Motivated by the $\Gamma$-convergence, the TV-$H^{-1}$ method is proposed. The inpainted image $u$ of $f\in L^2(\Omega)$ shall evolve via
\begin{align*}
u_t=\Delta p+\lambda\chi_{\Omega\backslash D}(f-u), ~~~p\in\partial \text{TV}(u)
\end{align*}
where $\partial$TV(u) denotes the subdifferential of 
\begin{align}
TV(u)=\left\{
                \begin{array}{ll}
                  |Du|(\Omega), ~~ \text{if $|u(x)|\leq1$ in $\Omega$}  \\
                  +\infty, ~~ otherwise
                \end{array}
              \right.
\end{align}
By using the existence of stationary solution in the Cahn-Hilliard case, it can be proved that 
\begin{align}
\Delta p+\lambda\chi_{\Omega\backslash D}(f-u)=0, ~~~p\in\partial TV(u)
\end{align}
admits a solution $u\in BV(\Omega)$

\section{Cahn-Hilliard equation on image inpainting}
Bertozzi et al. introduced the fourth order Cahn-Hilliard inpainting approach for binary, i.e., black and white, images~\cite{bertozzi}. This model is based on scalar smooth Cahn-Hilliard equation. The binary Cahn-Hilliard inpainting model has been generalized to gray value images~\cite{bosch}. This model is based on vector valued Cahn-Hilliard equation. (RGB for image of three channels, which corresponds three different chemical in fluid dynamics)

\subsection{Binary images}
Let f(x,y) be a given image in a domain $\Omega$, and suppose that $D\subset\Omega$ is the inpainting domain. Let $u(x,y,t)$ evolve in time to become a fully inpainted version of $f(x,y)$ under the equation:
\begin{align}
u_t=-\Delta(\varepsilon\Delta u-\frac{1}{\varepsilon}W'(u))+\lambda(x,y)(f-u)
\label{eq:cahn}
\end{align}
where
\begin{align*}
\lambda(x,y)=\left\{
                \begin{array}{ll}
                  0, ~~~\text{if } (x,y)\in D  \\
                  \lambda_0, ~~~\text{if } (x,y)\in \Omega \backslash D
                \end{array}
              \right.
\end{align*}
The function $W(u)$ is a nonlinear potential with wells corresponding to values of $u$ that are taken on by most of the grayscale values. In the binary case, $W$ should have wells at the values $u=0$ and $u=1$. We use the function $W(u)=u^2(u-1)^2$. We use convexity splitting.

The modified Cahn-Hilliard equation is not strictly a gradient flow. The original Cahn-Hilliard equation is indeed a gradient flow using an $H^{-1}$ norm for the energy.
\begin{align*}
E_1=\int_{\Omega}\frac{\varepsilon}{2}|\nabla  u|^2+\frac{1}{\varepsilon} W(u) dx
\end{align*}
The fidelity term of eq.~\ref{eq:cahn} can be derived from a gradient flow under an $L^2$ norm for the energy.
\begin{align}
E_2=\lambda_0\int_{\Omega \backslash D}(f-u)^2 dx
\end{align}

We can split $E_1$ as
\begin{align}
E_1=E_{11}-E_{12}
\end{align}
where
\begin{align*}
E_{11}&=\int_{\Omega} \frac{\varepsilon}{2}|\nabla u|^2+\frac{C_1}{2}|u|^2 dx \\
E_{12}&=\int_{\Omega}-\frac{1}{\varepsilon} W(u)+\frac{C_1}{2}|u|^2 dx
\end{align*}

A possible splitting for $E_2$ is
\begin{align*}
E_2=E_{21}-E_{22}
\end{align*}
where
\begin{align*}
E_{21}&=\int_{\Omega\backslash D} \frac{C_2}{2}|u|^2 dx \\
E_{22}&=\int_{\Omega\backslash D}-\lambda_0(f-u)^2+\frac{C_2}{2}|u|^2 dx
\end{align*}

For this splittings, the resulting time-stepping scheme is:
\begin{align*}
\frac{u^{n+1}-u^{n}}{\Delta t}=-\nabla_{H^{-1}}(E_{11}^{n+1}-E_{12}^n)-\nabla_{L^2}(E_{21}^{n+1}-E_{22}^{n})
\end{align*}
where $\nabla_{H^{-1}}$ and $\nabla_{L^2}$ represent gradient descent with respect to the $H^{-1}$ inner product, and $L^2$ inner product, respectively. This translates to a numerical scheme of the form
\begin{align*}
&\frac{u^{n+1}(x,y)-u^{n}(x,y)}{\partial t}+\varepsilon \Delta^2 u^{n+1}(x,y)-C_1\Delta u^{n+1}(x,y)+C_2 u^{n+1}(x,y) \\
=&\Delta(\frac{1}{\varepsilon}W'(u^{n}(x,y)))+\lambda(x,y)(f(x,y)-u^{n}(x,y))-C_1\Delta u^n(x,y)+C_2u^n(x,y)
\end{align*}

The constant $C_1$ and $C_2$ are positive, and need to be chosen large enough so that the energies $E_{11}, E_{12}, E_{21}$ and $E_{22}$ are convex. $C_1$ should be comparable to $\frac{1}{\varepsilon}$, while $C_2$ should be comparable to $\lambda_0$. We can use the spectral method and finite element method to solve this numerical PDE, to obtain an approximated solution of $u$.

\subsection{Color images}
Let $f$ be the given gray value image, which is defined on the image domain $\Omega\in\mathbb{R}^l$, with $l\in\{2,3\}$. Let $Y$ be the number of gray values which form the image. These $Y$ gray values are collected in the vector $\bold{g}=[g_1,g_2,\cdots,g_Y]^T \in\mathbb{R}^Y$, and $2\leq Y\leq 256$. The target is to reconstruct the image $f$ in the inpainting region $D$. Denote the reconstructed image by $f_r$. Let $T>0$ be a fixed time. Introducing a vector-valued phase variables $\bold{u}=[u_1,\cdots,u_Y]^T: \Omega\times(0,T)\rightarrow \mathbb{R}^Y$. $u_i$ describes the concentration of gray value $g_i$ for $i=1, \cdots, Y$. If $u_i(\bold{x},t)\approx 1$, then only gray value $g_i$ is present at point $\bold{x}$ at time $t$. $u_i(\bold{x},t)\approx 0$ means gray value $g_i$ is absent at point $\bold{x}$ at time $t$. Values of $u_i$ between $0$ and $1$ represent mixed regions. We initialize $u_i$ with $u_i(\bold{x},0)=f_i(\bold{x})$. The evolution of the reconstructed image $f_r$ is obtained from the components $u_i$ via,
\begin{align*}
f_r=\sum\limits_{i=1}^{Y}g_i u_i
\end{align*}
We have $\sum\limits_{i=1}^Y u_i=1$. The final reconstructed image $f_r$ of $f$ is $f_r(\bold{x},T)$.

The Cahn-Hilliard inpainting model is based on the Ginzburg-Landau energy $E_1$.
\begin{align*}
E_1(\bold{u})=\int_{\Omega}\frac{\varepsilon}{2}\sum\limits_{i=1}^{Y}|\nabla u_i|^2+\frac{1}{\varepsilon}W(\bold{u})d\bold{x}
\end{align*}

According to \cite{bosch}, the smooth gray value Cahn-Hilliard equation is
\begin{align*}
\left\{
                \begin{array}{lll}
&\partial_t u_i=\Delta w_i+\lambda(f_i-u_i) \\
&w_i=-\varepsilon \Delta u_i+\frac{1}{\varepsilon}W'(u_i)-\frac{1}{\varepsilon N}\sum\limits_{j=1}^N W'(u_j) \\
&\nabla u_i\cdot \bold{n}=\nabla w_i\cdot\bold{n}=0~~~\text{on~~}\partial\Omega
\end{array}
              \right.
\end{align*}
where $W(u)=u^2(u-1)^2$.

\section{Randomness, Stochastic Galerkin method}
When there is noise in the image, the initial condition of the Cahn-Hilliard equation will have uncertainty. We therefore introduce random variable $Z$ to the modified Cahn-Hilliard equation.

\subsection{Polynomial approximation}
Let $\mathbb{P}_n$ be the linear space of polynomials of degree at most $n$:
\begin{align*}
\mathbb{P}_n=\text{span}\{x^{k}: k=0,1,\cdots,n \}.
\end{align*}

\textbf{Weierstrass approximation theory}:
Let $I$ be a bounded interval and let $f\in C^{0}(\bar{I})$. Then, for any $\epsilon > 0$, we can find $n\in\mathbb{N}$ and $p\in\mathbb{P}_n$ such that
\begin{align*}
|f(x)-p(x)|<\epsilon,~~~~~~\forall x\in\bar{I}
\end{align*}

In other words, we would like to study the existence of $\phi_n(f)\in\mathbb{P}_n$ such that
\begin{align}
||f-\phi_n(f)||_{\infty}=\inf_{\substack{\psi\in\mathbb{P}_n}} ||f-\psi||_{\infty}
\label{eq:norm}
\end{align}

The $n$th-degree polynomial $\phi_n(f)$ is called the polynomial of best uniform approximation of $f$ in $\bar{I}$.

Another approximation problem can be formulated in terms of norms other than the infinity norm used in eq.~(\ref{eq:norm}). For a positive weight function $w(x)$, $x\in I$, the weighted $L^2$ space by:
\begin{align*}
L_w^2(I)=\left\{v: I\rightarrow\mathbb{R}| \int_I v^2(x)w(x)dx<\infty  \right\}
\end{align*}

with the inner product
\begin{align*}
(u,v)_{L_w^2(I)}=\int_I u(x)v(x)w(x)dx,~~~~~~\forall u,v\in L_w^2(I),
\end{align*}
and the norm
\begin{align*}
||u||_{L_w^2(I)}=\left( \int_I u^2(x)w(x)dx\right)^{1/2}.
\end{align*}

\subsubsection{Orthogonal Projection}
Let $N$ be a fixed nonnegative integer and let $\{\phi_k(x)\}_{k=0}^{N}\subset\mathbb{P}_N$ be orthogonal polynomials of degree at most $N$ with respect to the positive weight $w(x)$.
\begin{align*}
(\phi_m(x),\phi_n(x))_{L_w^2(I)}=||\phi_m||^2_{L_w^2(I)}\delta_{m,n},~~~~~0\leq m,n\leq N.
\end{align*}

The projection operator $P_N: L_w^2(I)\rightarrow\mathbb{P}_N$, for any function $f\in L_w^2(I)$
\begin{align*}
P_N f=\sum\limits_{k=0}^{N}\hat{f}_k\phi_k(x)
\end{align*}
where
\begin{align*}
\hat{f}_k=\frac{1}{||\phi_k||^2_{L_w^2}}(f,\phi_k)_{L_w^2},~~~~~~~0\leq k\leq N
\end{align*}

Obviously, $P_N f\in\mathbb{P}_N$. It is called the orthogonal projection of $f$ onto $\mathbb{P}_N$ via the inner product $(\cdot,\cdot)_{L_w^2}$, and $\{\hat{f}_k\}$ are the generalized Fourier coefficients.

The following trivial facts hold:
\begin{align*}
P_N f&=f,~~~~~\forall f\in \mathbb{P}_N \\
P_N\phi_k&=0,~~~~~\forall k>N.
\end{align*}

\textbf{Theorem}: For any $f\in L_w^2(I)$ and any $N\in\mathbb{N}_0$, $P_N f$ is the best approximation in the weighted $L^2$ norm in the sense that
\begin{align*}
||f-P_N f||_{L_w^2}=\inf_{\substack{\psi\in\mathbb{P}_n}} ||f-\psi||_{L_w^2}
\end{align*}

Proof. Any polynomial $\psi\in\mathbb{P}_N$ can be written in the form $\psi=\sum\limits_{k=0}^N c_k\phi_k$ for some real coefficients $c_k$, $0\leq k\leq N$. Minimizing $||f-\psi||_{L_w^2}^2$, whose derivatives are:
\begin{align*}
\frac{\partial}{\partial c_j}||f-\psi||_{L_w^2}^2&= \frac{\partial}{\partial c_j}\left( ||f||_{L_w^2}^2-2\sum\limits_{k=0}^{N}c_k(f,\phi_k)_{L_w}^2+\sum\limits_{k=0}^{N}c_k^2||\phi_k||_{L_w^2}^2 \right) \\
&=-2(f,\phi_j)_{L_w^2}+2c_j||\phi_j||_{L_w^2}^2, ~~~~~0\leq j\leq N.
\end{align*}

By setting the derivatives to zero, the unique minimum is attained when $c_j=\hat{f}_j,~~0\leq j\leq N$, where $\hat{f}_j$ are the Fourier coefficients of $f$. This completes the proof.

\subsubsection{Spectral Convergence}
The convergence of the orthogonal projection can be stated as follows:\\
\textbf{Theorem}: For any $f\in L_w^2(I)$,
\begin{align*}
\lim_{N\rightarrow\infty}||f-P_N f||_{L_w^2}=0.
\end{align*}

The rate of convergence depends on the regularity of $f$ and the type of orthogonal polynomial $\{\phi_k\}$. Defined a weighted Sobolev space $H_m^{k}(I)$, for $k=0,1,2,\cdots$, by
\begin{align*}
H_w^k(I)=\left\{v: I\rightarrow \mathbb{R} | \frac{d^m v}{dx^m}\in L_w^2(I), 0\leq m\leq k\right\}
\end{align*}
equipped an inner product
\begin{align*}
(u,v)_{H_w^k}=\sum\limits_{m=0}^k\left(\frac{d^mv}{dx^m},\frac{d^mv}{dx^m} \right)_{L_w^2}
\end{align*}
and a norm $||u||_{H_w^k}=(u,u)_{H_w^k}^{1/2}$.

Consider the case of $\bar{I}=[-1,1]$ with weight function $w(x)=1$ and Legendre polynomials $\{P_n(x)\}$. The orthogonal projection for any $f(x)\in L_w^2(I)$ is
\begin{align*}
P_N f(x)&=\sum\limits_{k=0}^{N}\hat{f}_{k}P_k(x) \\
\hat{f}_k&=\frac{1}{||P_k||_{L_w^2}^2}(f,P_k)_{L_w^2}.
\end{align*}

The following result holds.

\textbf{Theorem.} For any $f(x)\in H_w^p[-1,1], p\geq 0$, there exists a constant C, independent of $N$, such that
\begin{align*}
||f-P_N f||_{L_w^2[-1,1]}\leq CN^{-p} ||f||_{H_w^p[-1,1]}
\end{align*}

Since the Legendre polynomial satisfy
\begin{align*}
Q(P_k)=\lambda_k P_k
\end{align*}
where
\begin{align*}
Q=\frac{d}{dx}\left((1-x^2)\frac{d}{dx} \right)=(1-x^2)\frac{d^2}{dx^2}-2x\frac{d}{dx}
\end{align*}
and $\lambda_k=-k(k+1)$. We then have
\begin{align*}
(f,P_k)_{L_w}^2&=\frac{1}{\lambda_k}\int_{-1}^{1}Q[P_k]f(x)dx \notag \\
&=\frac{1}{\lambda_k}\int_{-1}^{1}((1-x^2)P_k'' f-2xP_k'f)dx \notag \\
&=-\frac{1}{\lambda_k}\int_{-1}^{1}[((1-x^2)f)'P_k' +2xP_k'f)]dx \notag \\
&=-\frac{1}{\lambda_k}\int_{-1}^{1}(1-x^2)f'P_k' dx \notag \\
&=\frac{1}{\lambda_k}\int_{-1}^{1}((1-x^2)f')'P_k dx
\end{align*}
where the third and the last equality use the rule of integration by parts. This implies
\begin{align*}
(f,P_k)_{L_m^2}=\frac{1}{\lambda_k}(Q[f],P_k)_{L_w^2}
\end{align*}

By applying the procedure repeatedly for $m$ times, we have
\begin{align*}
(f,P_k)_{L_m^2}=\frac{1}{\lambda_k^m}(Q^m[f],P_k)_{L_w^2}
\end{align*}

The projection error can be estimated as
\begin{align*}
||f-P_N f||_{L_w^2}^2&=\sum\limits_{k=N+1}^{\infty}\hat{f}_k^2||P_k||_{L_w^2}^2 \notag\\
&=\sum\limits_{k=N+1}^{\infty}\frac{1}{||P_k||^2_{L_w^2}}(f,P_k)_{L_w^2}^2 \notag \\
&=\sum\limits_{k=N+1}^{\infty}\frac{1}{\lambda_k^{2m}||P_k||^2_{L_w^2}}(Q^m[f],P_k)_{L_w^2}^2 \notag \\
&\leq \lambda_N^{-2m}\sum\limits_{k=0}^{\infty}\frac{1}{||P_k||^2_{L_w^2}}(Q^m[f],P_k)_{L_w^2}^2 \notag \\
&\leq  N^{-4m}||Q^m[f]||_{L_w^2}^2\leq CN^{-4m}||f||_{H_w^{2m}}^2
\end{align*}

\subsection{Generalized Polynomial Chaos (gPC)}
In the gPC expansion, one approximates the solution of a stochastic problem via an orthogonal polynomial series.

\subsubsection{Multiple random variables}
Let $Z=(Z_1,\cdots, Z_d)$ be a random vector with mutually independent components and distribution $F_Z(z_1,\cdots, z_d)=P(Z_1\leq z_1,\cdots, Z_d\leq z_d)$. For each $i=1,\cdots, d$, let $F_{Z_i}(z_i)=P(Z_i\leq z_i)$ be the marginal distribution of $Z_i$, whose support is $I_{Z_i}$. Mutual independence among all $Z_i$ implies that $F_Z(z)=\Pi_{i=1}^d F_{Z_i}(z_i)$ and $I_Z=I_{Z_1}\times\cdots\times I_{Z_d}$. Also, let $\{\phi_k(Z_i)\}_{k=0}^N\in\mathbb{P}_N(Z_i)$ be the univariate gPC basis functions in $Z_i$ of degree up to $N$. That is,
\begin{align*}
\mathbb{E}[\phi_m(Z_i)\phi_n(Z_i)]=\int\phi_m(z)\phi_n(z)dF_{Z_i}(z)=\delta_{mn}\gamma_{m},~~~~~0\leq m,n\leq N
\end{align*}

Let $\bold{i}=(i_1,\cdots, i_d)\in\mathbb{N}_0^d$ be a multi-index with $\bold{i}=i_1+\cdots+i_d$. Then the $d$-variate $N$th-degree gPC basis functions are the products of the univariate gPC polynomials of total degree less than or equal to $N$ (tensor product of basis functions):
\begin{align*}
\Phi_{\bold{i}}(Z)=\phi_{i_1}(Z_1)\cdots\phi_{i_d}(Z_d), ~~~~~~~0\leq|\bold{i}|\leq N.
\end{align*}

It follows immediately that
\begin{align*}
\mathbb{E}[\Phi_{\bold{i}}(Z)\Phi_{\bold{j}}(Z)]=\int\Phi_{\bold{i}}(z)\Phi_{\bold{j}}(z)dF_Z(z)=\gamma_{\bold{i}}\delta_{\bold{ij}}
\end{align*}
where $\gamma_{\bold{i}}=\mathbb{E}[\Phi_{\bold{i}}^2]=\gamma_{i_1}\cdots\gamma_{i_d}$ are the normalization factors and $\delta_{\bold{ij}}=\delta_{i_1 j_1}\cdots\delta_{i_d j_d}$ is the $d$-variate Kronecker delta function. It is obvious that the span of the polynomials is $\mathbb{P}_N^d$, the linear space of all polynomials of degree at most $N$ in $d$ variables.
\begin{align*}
\mathbb{P}_N^d(Z)=\left\{p: I_Z\rightarrow\mathbb{R}~|~p(Z)=\sum\limits_{|\bold{i}|\leq N}c_{\bold{i}}\Phi_{\bold{i}}(Z) \right\}
\end{align*}
whose dimension is
\begin{align*}
\dim\mathbb{P}_N^d= \binom {N+d} N ,
\end{align*}

The $d$-variate gPC projection follows the univariate projection in a direct manner. Let $L_{dFz}^2(I_Z)$ be the space of all mean-square integrable functions of $Z$ with respect to the measure $dF_Z$, that is:
\begin{align*}
L_{dF_Z}=\left\{f:I_Z\rightarrow\mathbb{R}~|~ \mathbb{E}[f^2(Z)]=\int_{I_Z}f^2(z)dF_Z(z)<\infty \right\}
\end{align*}

Then for $f\in L_{dF_Z}^2$, its $N$th-degree \textbf{gPC orthogonal projection} is defined as
\begin{align}
P_N f=\sum\limits_{|\bold{i}|\leq N}\hat{f}_{\bold{i}}\Phi_{\bold{i}}(Z)
\end{align}

where
\begin{align*}
\hat{f}_{\bold{i}}=\frac{1}{\gamma_{\bold{i}}}\mathbb{E}[f\Phi_{\bold{i}}]=\frac{1}{\gamma_{\bold{i}}}\int f(z)\Phi_{\bold{i}}(z)dF_Z(z),~~~~~\forall |\bold{i}|\leq N.
\end{align*}
The classical approximation theory can be readily applied to obtain
\begin{align*}
||f-P_N f||_{L_{dF_Z}^2}\rightarrow 0, ~~~~~N\rightarrow\infty,
\end{align*}
and
\begin{align*}
||f-P_N f||_{L_{dF_Z}^2}=\inf_{\substack{g\in \mathbb{P}_N^d}} ||f-g||_{L_{dF_Z}^2}
\end{align*}

The correspondence between the type of generalized Polynomial Chaos and their underlying random variables (distribution of $Z$) is shown in Table~\ref{table:gpc}.
\begin{table}[ht]
\centering
\caption{Generalized Polynomial Chaos (gPC) and their underlying random variables}
\label{table:gpc}
\begin{tabular}{llll}
\hline \hline
 & Distribution of $Z$ & gPC basis polynomials  & Support   \\
 \hline
 Continuous & Gaussian  & Hermite  &  $(-\infty,\infty)$  \\
 &  Gamma & Laguerre  & $[0,\infty)$   \\
 & Beta & Jacobi & $[a,b]$  \\
 & Uniform & Legendre & $[a,b]$ \\ [1ex]
 Discrete & Poisson & Charlier & $\{0,1,2,\cdots\}$ \\
 & Binomial & Krawtchouk &$\{0,1,\cdots, N\}$ \\
 & Negative binomial & Meixner & $\{0,1,2,\cdots\}$ \\
 & Hypergeometric & Hahn & $\{0,1,\cdots, N\}$ \\ [1ex]
 \hline
\end{tabular}
\end{table}

\subsection{Parametric Cahn-Hilliard equation}
The stochastic Cahn-Hilliard equation can be formulated as follows:
\begin{align}
\left\{
                \begin{array}{llll}
&\partial_t u(x,y,t,Z)=\Delta w(x,y,t,Z)+\lambda(x,y)(f-u(x,y,t,Z)),~~~\text{on~ } \Omega\times(0,T]\times \mathbb{R}^d    \\
&w(x,y,t,Z)=-\varepsilon \Delta u(x,y,t, Z)+\frac{1}{\varepsilon}W'(u),~~~\text{on~ } \Omega\times(0,T]\times \mathbb{R}^d \\
& \nabla u\cdot n =\nabla w\cdot n=0, ~~~\text{on~ } \partial\Omega \\
&u(x,y,0,Z)=u^{0}(x,y,Z)=u^{0}(x,y)+Z
\end{array}
              \right.
\label{eq:cahn_random}
\end{align}
where $Z=(Z_1,Z_2,\cdots, Z_d)\in \mathbb{R}^d,~~d\geq 1$ are a set of mutually independent random variables characterizing the random inputs to the governing equation. We are trying to compute the statistics of the solution when uncertainty is involved in the system of equations.  The generalized polynomial chaos expansion is a representation of stochastic processes by polynomial functionals of random variables:
\begin{align*}
u(x,y,t,Z)=\sum\limits_{i=0}^{\infty}u_i(x,y,t)\Phi_i(Z)
\end{align*}

The finite-term expansion takes the form:
\begin{align}
u_N(x,y,t,Z)=\sum\limits_{i=0}^{N}u_i(x,y,t)\Phi_i(Z)
\label{eq:expansion}
\end{align}
where $N$ is the highest order of the expansion. Substituting eq.~(\ref{eq:expansion}) into eq.~(\ref{eq:cahn_random}), we obtain
\begin{align}
\left\{
                \begin{array}{llll}
&\sum\limits_{i=0}^{N}\frac{\partial u_i(x,y,t)}{\partial t}\Phi_i(Z)=\Delta\left(\sum\limits_{i=0}^{N} w_i(x,y,t)\Phi_i(Z)\right)+\lambda(x,y) \left(f-\sum\limits_{i=0}^{N}u_i(x,y,t)\Phi_i(Z)\right) \\
&\sum\limits_{i=0}^{N} w_i(x,y,t)\Phi_i(Z)=-\varepsilon\Delta\left( \sum\limits_{i=0}^{N}u_i(x,y,t)\Phi_i(Z)\right)+\frac{1}{\varepsilon}\sum\limits_{i=0}^{N}W'(u)\Phi(Z) \\
&\nabla\left(\sum\limits_{i=0}^{N}u_i\Phi(Z)\right)\cdot n=\nabla\left(\sum\limits_{i=0}^{N}w_i\Phi(Z)\right)\cdot n=0
\end{array}
              \right.
\label{eq:cahn_expan}
\end{align}

Projecting the above equation onto the bases spanned by $\{\Phi_j\}_{j=0}^N$, and using the orthogonality of the bases:
\begin{align*}
<\Phi_i\Phi_j>=<\Phi_i^2>\delta_{ij}
\label{eq:ortho}
\end{align*}
where $\delta_{ij}$ is the Kroncker delta and $<\cdot,\cdot>$ denotes the ensemble average which is the inner product in the Hilbert space of the variable $Z$, we obtain, according to~\cite{d_xiu}:
\begin{align}
\left\{
                \begin{array}{lll}
&\mathbb{E}(\frac{\partial u_N}{\partial t}\Phi_{\bold{j}})=\mathbb{E}(\Delta w_N \Phi_{\bold{j}})+\lambda(\mathbb{E}(f-u_N)\Phi_j) \\
&\mathbb{E}(w_N\Phi_\bold{j})=-\varepsilon\mathbb(\Delta u_N\Phi_{\bold{j}})+\frac{1}{\varepsilon}\mathbb{E}((4u_N^3-6u_N^2+2u_N)\Phi_j) \\
&\mathbb{E}(\nabla u_N\Phi_{\bold{j}}\cdot n)=\mathbb{E}(\nabla w_N\Phi_{\bold{j}}\cdot n)=0
\end{array}
              \right.
\end{align}
therefore,
\begin{align}
\left\{
                \begin{array}{lll}
&\frac{\partial u_j}{\partial t}=\Delta w_j+\lambda(f-u_j) \\
&w_j=-\varepsilon\Delta u_j+\frac{1}{\varepsilon} \left(4\mathbb{E}(u^3\Phi_{\bold{j}})-6\mathbb{E}(u^2\Phi_{\bold{j}})+2\mathbb{E}(u\Phi_{\bold{j}})\right) \\
&\nabla u_j\cdot n=\nabla w_j\cdot n
\end{array}
              \right.
\label{eq:ch_gpc2}
\end{align}

The second equation of eq.~(\ref{eq:ch_gpc2}) will lead us to:
\begin{align}
w_j=-\varepsilon\Delta u_j+\frac{1}{\varepsilon}\left(4\left(\sum\limits_{i=0}^{N}\sum\limits_{p=0}^{N}\sum\limits_{q=0}^{N} u_i u_p u_q e_{ipqj}\right)/\gamma_1-6\left(\sum\limits_{i=0}^{N}\sum\limits_{p=0}^{N}u_i u_p e_{ipj}\right)/\gamma_2+2u_j\right)
\end{align}
where $e_{ipqj}=\mathbb{E}(\Phi_i\Phi_p\Phi_q\Phi_j)$, $e_{ipj}=\mathbb{E}(\Phi_i\Phi_p\Phi_j)$, $\gamma_1=\gamma_2=\mathbb{E}(\phi_j^2)$.

When the noise is Gaussian, and is related to just one random variable, to simplify the problem, let $Z\sim\mathcal{N}(0,\sigma^2)$,  the initial condition can be written as:
\begin{align*}
u(x,y,0,Z)=u^0(x,y)+Z
\end{align*}
Under this condition, the generalized polynomial chaos expansion for the initial condition takes the form
\begin{align}
u_0^0=u^0(x,y),~~~u_1^0=1, ~~~u_k^0=0, \text{ for } k\geq 2
\end{align}
Since we have
\begin{align*}
\left\{
                \begin{array}{lllll}
&e_{0000}=\int\rho(z)dz=1; \\
&e_{0001}=e_{1000}=e_{0100}=e_{0010}=e_{100}=e_{010}=e_{001}=\int z\rho(z)dz=0; \\
&e_{1001}=e_{0101}=e_{0011}=e_{1100}=e_{0110}=e_{1100}=e_{110}=e_{011}=\int z^2\rho(z)dz=C_1 \\
&e_{1110}=e_{0111}=e_{111}=\int z^3\rho(z)dz=0 \\
&e_{1111}=\int z^4\rho(z)dz=C_2
\end{array}
              \right.
\end{align*}
where $\rho(z)$ is the distribution of random variable $Z$,
therefore, eq.~(\ref{eq:ch_gpc2}) reduces to
\begin{align}
&\frac{\partial u_0}{\partial t}=-\Delta\Biggl(\varepsilon\Delta u_0- \frac{1}{\varepsilon} \biggl( 4(u_0^3+3u_0u_1^2\cdot C_1)-6(u_0^2+u_1^2\cdot C_1)+2u_0 \biggr)  \Biggr) +\lambda(f-u_0) \\
&\frac{\partial u_1}{\partial t}=-\Delta\Biggl(\varepsilon\Delta u_1- \frac{1}{\varepsilon} \biggl( 4(u_1^3\cdot C_2+3u_1u_0^2\cdot C_1)-6(2u_0u_1\cdot C_1)+2u_1 \biggr)  \Biggr) -\lambda u_1
\end{align}
with initial condition $u_0^0=u^0(x,y)$ and $u_1^0=1$.

\section{Framework of using wavelet}
Take the case when the noise is of uniform distribution in $(-1, 1)$ as an example. Let $z$ be the random variable. The pdf is $\rho(z)=\frac{1}{2}$ and is a constant. The orthogonality of bases defines the Legendre orthogonal polynomials:
\begin{align}
\int_{-1}^{1} P_n(z)P_m(z)dx=\frac{2}{2n+1}\delta_{nm},
\end{align}
and
\begin{align}
P_{n+1}=\frac{2n+1}{n+1}
\end{align}
with,
\begin{align*}
P_0(z)=1,~~~P_1(z)=z,~~~P_2(z)=\frac{3}{2}z^2-\frac{1}{2},~\cdots
\end{align*} 
 
\begin{figure}[!ht]
% \hfill
\center
\includegraphics[width=.48\linewidth,height=3.5cm]{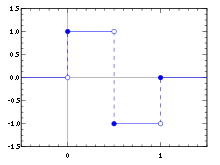}
\caption{Haar Wavelet, $\beta(x)$}
\label{fig:haar}
\end{figure}

 The orthonormal polynomial system $\Phi=\{\phi_0, \phi_1, \phi_2, \cdots\}=(1,~z,~\frac{3}{2}z^2-\frac{1}{2},~\cdots)$. In wavelet system
 \begin{align*}
 \Phi&=\{\varphi_{0,k}: k\in\mathbb{Z}\}\bigcup\{\psi_{j,k}: k\in\mathbb{Z}\}_{j=0}^{\infty} \\
 &=\{ \varphi_{0,0}, \varphi_{0,1}, \varphi_{0,-1}, \varphi_{0,2}, \varphi_{0,-2}, \cdots \} \bigcup \{ \psi_{0,0}, \psi_{0,1}, \psi_{0,-1}, \cdots \}  \bigcup \{ \psi_{j,0}, \psi_{j,1}, \psi_{j,-1}, \cdots \}_{j=1}^{\infty} 
 \end{align*}
where,
\begin{align*}
\varphi_{0,k}(x)&=\chi_{[0,1)}(x-k) \\
\psi_{j,k}(x)&=c_j\beta(2^j-k)
\end{align*}
For Haar wavelet,
\begin{align}
\beta(x)=\left\{
                \begin{array}{lll}
&1,~~0\leq x< \frac{1}{2} \\
&-1,~~\frac{1}{2}\leq x<1 \\
& 0, ~~\text{otherwise}
\end{array}
              \right.
\end{align} 
and,
\begin{align}
\beta(2x)=\left\{
                \begin{array}{lll}
&1,~~0\leq x< \frac{1}{4} \\
&-1,~~\frac{1}{4}\leq x<\frac{1}{2} \\
& 0, ~~\text{otherwise}
\end{array}
              \right. 
              \\
\beta(2x-1)=\beta(2(x-\frac{1}{2}))=
\left\{
\begin{array}{lll}
&1,~~\frac{1}{2}\leq x< \frac{3}{4} \\
&-1,~~\frac{3}{4}\leq x<1 \\
& 0, ~~\text{otherwise}
\end{array}
              \right.
\end{align}
$\Phi$ constitutes an orthonormal bases in $L^2(\mathbb{R})$. Figure~\ref{fig:haar} shows a plot of standard Haar wavelet $\beta(x)$.

The connection of polynomial bases and Haar wavelet is shown in Table~\ref{table:poly_wavelet}. $\varphi_{0,k}$ is a coarse scale representation, while $\psi_{j,k}$ is a finer scale representation. In high dimension, tensor product will be used. 

\begin{table}[ht]
\centering
\caption{Connection of Polynomial Chaos (PC) and the wavelet basis}
\label{table:poly_wavelet}
\begin{tabular}{lll}
\hline \hline
Legendre polynomial & wavelet basis  & wavelet expression\\
 with variable $z$ & ~  \\
 \hline
 \centering
  1  &  $\varphi_{0,k}$ & $\chi_{[0,1)}(z-k)$   \\
  $z$ &    $\psi_{0,k}$  & $c_0\beta(z-k)$\\
 $\frac{3}{2}z^2$ & $\psi_{1,k}$ &$c_1\beta(2z-k)$ \\
 $\cdots$ & $\cdots$ & $\cdots$ \\
 [1ex]
 \hline
\end{tabular}
\end{table}

Projecting a variable $Z$ onto the wavelet basis, we have
\begin{align}
Z=\sum\limits_{k\in\mathbb{Z}}<Z,\varphi_{0,k}>\varphi_{0,k}+\sum\limits_{j=0}^{\infty}\sum\limits_{k\in\mathbb{Z}}<Z,\psi_{j,k}>\psi_{j,k}
\end{align}

Considering the wavelet with $n$ vanishing moment, that is 
\begin{align}
\int_{-\infty}^{\infty}z^k\psi(z)dz=0,~~(0\leq k<n),
\end{align}
when $n=2$, we have $\int_{-\infty}^{\infty}z\psi(z)dz=0$, and $\int_{-\infty}^{\infty}z^2\psi(z)dz=0$. (We can express it in the form of Daubechies wavelet, although there is no explicit expression for the wavelet terms $\psi$. Haar wavelet has an explicit expression.) Referring to the polynomial case, we can use wavelet to express the equation terms, for example: $e_{ipqj}=\mathbb{E}(\Phi_i\Phi_p\Phi_q\Phi_j)=\int\varphi_{0,k}(z)\psi_{j,k}(z)\varphi_{0,k}(z)\varphi_{0,k}(z)\rho(z)dz$. When the parameters of uncertainty is high dimension, the computation is still difficult.

\section{Perturbation method}
When $Z\in(0,\delta)$ is a random variable, and $\delta<<1$. The initial condition is uncertain, and can be written as
\begin{align*}
u^0(x,y,Z)=u(x,y,0,Z)=u^0(x,y)+Z
\end{align*}
In the perturbative approach, the stochastic quantities are expanded via a Taylor series around the mean value of the random  inputs,
\begin{align}
u^0(x,y,Z)=u(x,y,0,Z)=u_0^0(x,y)+Z u_1^0(x,y)+\cdots
\label{eq:expansion}
\end{align}
where 
\begin{align*}
u_k(x,y,t)=k!\frac{\partial^k u(x,y,t,Z)}{\partial^k Z}
\end{align*}
the expansion is at the point of the mean of the distribution. Substituting expansion~(\ref{eq:expansion}) into the Cahn-Hilliard equation~(\ref{eq:cahn}), and equating the terms of different orders, under the assumption that $O(1)>>O(Z)>>O(Z^2)>>\cdots$,  and the condition that $u_0^0=u_0$, $u_1^0=1$, and $u_k^0=0$ for $k\geq 2$, we have
\begin{align}
&O(Z^0):~~\frac{\partial u_0}{\partial t}=-\Delta(\varepsilon \Delta u_0-\frac{1}{\varepsilon}(4u_0^3-6u_0^2+2u_0)) +\lambda(f-u_0) \\
&O(Z^1):~~\frac{\partial u_1}{\partial t}=-\Delta(\varepsilon\Delta u_1-\frac{1}{\varepsilon}(4\cdot(3u_0^2u_1)-6(2u_0u_1)+2u_1))-\lambda u_1 
\end{align} 

\section{Experiments}
\begin{figure}[!ht]
% \hfill
\center
\subfigure[Masked image]{\includegraphics[width=.48\linewidth,height=3cm]{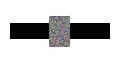}}
% \hfill
\subfigure[Output image]{\includegraphics[width=.48\linewidth,height=3cm]{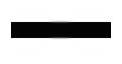}}
\caption{Binary image inpainting}
\label{fig:img_in}
\end{figure}

As shown in~\cite{bertozzi}, the effect of using binary Cahn-Hilliard equation for image inpainting is shown in Figure~\ref{fig:img_in}.

For Gaussian noise case, that is $Z\sim\mathcal{N}(0,1)$, with probability density function
\begin{align}
\rho(z)=\frac{1}{\sqrt{2\pi}}e^{-z^2/2},
\end{align}
applying the Stochastic Galerkin method for Cahn-Hilliard equation. According to Table~\ref{table:gpc}, and the orthogonality condition~\ref{eq:ortho}, the hermite orthogonal polynomial $\{H_m(Z)\}$:
\begin{align*}
H_0(Z)=1,~~~H_1(Z)=Z,~~~ H_2(Z)=Z^2-1,~~~ H_3(Z)=Z^3-3Z,~ \cdots
\end{align*}
is used as the basis function. 

\begin{figure}[!ht]
\center
\subfigure[Mean image of the final stage]{\includegraphics[width=.48\linewidth,height=3cm]{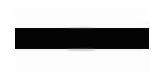}}
\subfigure[Mean image of the final stage]{\includegraphics[width=.48\linewidth,height=3cm]{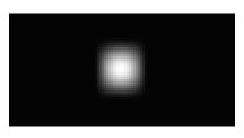}} \\
\caption{Image inpainting with stochastic Galerkin method}
\label{fig:img_ga}
\end{figure}

Applying the stochastic Galerkin method, we can obtain the mean and first order solution of $u$, as shown in Figure~\ref{fig:img_ga}.

\begin{figure}[!ht]
\center
\subfigure[Plot of leading order solution of the final stage]{\includegraphics[width=.48\linewidth,height=3cm]{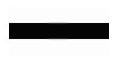}} 
\subfigure[First order solution of the final stage]{\includegraphics[width=.48\linewidth,height=3cm]{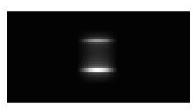}}
\caption{Image inpainting with perturbation expansion}
\label{fig:img_per}
\end{figure}

According to perturbation method, the mean solution and the first order solution are shown in Figure~\ref{fig:img_per}.

%\section{Summary}
%In this paper we developed a sparse wavelets based stochastic Galerkin method for the Cahn-Hilliard equation with uncertainty. The uncertainty comes from initial data. This method enables us to quantify the uncertainty from multi-dimensional random inputs, which is previously infeasible using the global gPC basis. We proved and numerically demonstrated the sparsity of the basis related coefficient, Sijk, which allows us to dramatically accelerate the computation of the collision operator under the Galerkin projection. Regularity of the solution of the Cahn-Hilliard equation in the random space and an accuracy result of the stochastic Galerkin method are proved.

%Many related problems are still open, for example, asymptotic-preserving schemes [10] for the Boltzmann equation with uncertainty, adaptive mesh techniques to capture discontinuities in the random space, quantification of nonlinear uncertainties on the collision kernel, etc
 
%\clearpage

\end{document}